%% file: main.tex
\pdfoutput=1
\documentclass[pdflatex,sn-mathphys-num]{sn-jnl}

\input{shared}

\hypersetup{
  pdftitle={Quadratic Envelope Search for Evaluation-Efficient Derivative-Free Optimization},
  pdfauthor={F. Demirci}
}

\begin{document}

\title[Quadratic Envelope Search]{Quadratic Envelope Search for
Evaluation-Efficient Derivative-Free Optimization}

\author*[1]{\fnm{Fatih} \sur{Demirci}}\email{f.demirci24@imperial.ac.uk}
\author[2]{\fnm{Lorenzo} \sur{Sabug Jr.}}\email{l.sabug21@imperial.ac.uk}
\author[1]{\fnm{Urban} \sur{Fasel}}\email{u.fasel@imperial.ac.uk}
\author[1,2]{\fnm{Eric C.} \sur{Kerrigan}}\email{e.kerrigan@imperial.ac.uk}

\affil*[1]{\orgdiv{Department of Aeronautics}, \orgname{Imperial College
London}, \orgaddress{\city{London}, \country{United Kingdom}}}
\affil[2]{\orgdiv{Department of Electrical and Electronic Engineering},
\orgname{Imperial College London}, \orgaddress{\city{London},
\country{United Kingdom}}}

\abstract{\unboldmath
We introduce Quadratic Upper--Lower Envelope Search (QULES), a derivative-free
method for unconstrained nonconvex optimization. A nonconvex subproblem seeks
to reduce the curvature gap between quadratic envelopes that bound sampled
objective values and share a center in a trust region around the incumbent.
QULES evaluates the objective at this center, where the upper envelope attains
its minimum, and compares the actual decrease with the envelope's prediction.
To represent curvature away from the coordinate axes, we introduce QULES-R, a
variant that adapts the envelopes' orthogonal frame using previously evaluated
points without needing additional objective evaluations.
We prove that the gradient norm tends to zero along the full sequence for
continuously differentiable objectives with a compact initial lower level set
and uniformly bounded upper-envelope curvature on radius reductions. With the
same curvature bound, Lipschitz gradients on sample segments and objective
values bounded below along the iterates, we bound the objective evaluations
needed to reach a gradient norm of at most $\varepsilon$. Including
initialization and geometry repair, the bounds are
$\mathcal O(n^2\varepsilon^{-2})$ for QULES with $\mathcal O(n)$ retained
samples and $\mathcal O(n^3\varepsilon^{-2})$ for QULES-R with
$\mathcal O(n^2)$ retained samples, when the remaining constants are
independent of $n$.
Within $50(n+1)$ evaluations per problem at a relative objective-gap tolerance
of $10^{-7}$, QULES-R reached 25/25 targets on rotated ill-conditioned problems
and 30/35 on the Mor\'e--Garbow--Hillstrom suite, compared with 13/25 and
26/35, respectively, for Powell's quadratic-interpolation method NEWUOA.}



\enlargethispage{8pt}
\maketitle

\input{sections/introduction}
\input{sections/problem-setting}
\input{sections/envelope-search-algorithm}
\input{sections/convergence-analysis}
\input{sections/numerical-experiments}

\input{sections/conclusion}

\appendix
\begin{figure}[!t]
\centering
\begingroup
\small
\captionof{table}{Dimensions of the 135 CUTEst problems. Standard starting points are used.}
\label{tab:cutest-dimensions}
\begin{tabular}{lrrrrrrrrrrr}
\toprule
Dimension $n$ & 2 & 3 & 4 & 5 & 6 & 7 & 8 & 9 & 10 & 11 & 12\\
\midrule
Number of problems & 49 & 29 & 11 & 5 & 7 & 8 & 16 & 2 & 6 & 1 & 1\\
\bottomrule
\end{tabular}
\par
\endgroup
\medskip
\includegraphics[width=\linewidth]{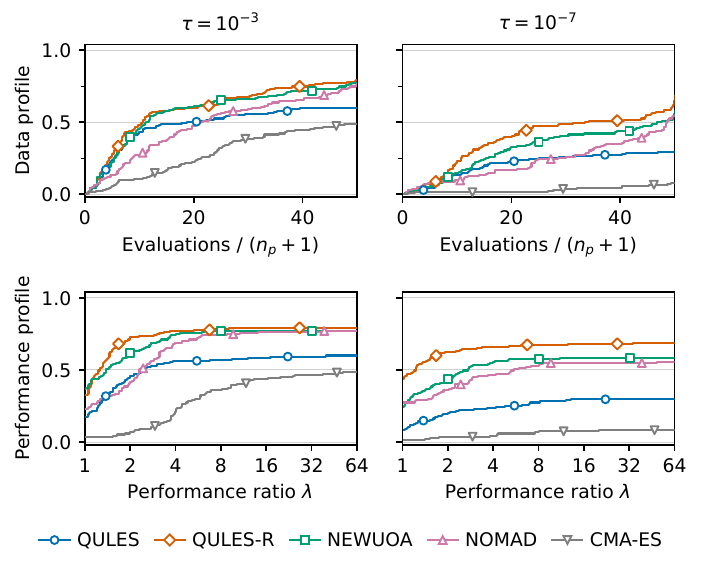}
\caption{135 CUTEst problems. Data profiles (top) and performance profiles
(bottom) at $\tau=10^{-3}$ (left) and $\tau=10^{-7}$ (right).}
\label{fig:cutest-profiles}
\end{figure}

\input{sections/frame-update}
\input{sections/cutest-appendix}


\section*{Statements and Declarations}

\noindent\textbf{Funding}
Fatih Demirci was supported by the Republic of T\"urkiye Ministry of National
Education. Lorenzo Sabug Jr. was supported by UKRI EPSRC grant EP/Z533816/1,
``Concurrent Learning and Control of Uncertain Large-Scale Phenomena''.

\noindent\textbf{Competing interests}
The authors have no competing interests to declare that are relevant to the
content of this article.

\noindent\textbf{Data and code availability}
The code, reproduction instructions and raw solver histories are available at
\url{https://doi.org/10.5281/zenodo.24572260}.

\input{references.bbl}
\end{document}

%% file: shared.tex
\usepackage{graphicx}
\usepackage{amsmath,amssymb,amsfonts}
\usepackage{amsthm}
\usepackage{mathtools}
\usepackage{bm}
\usepackage{float}
\usepackage{placeins}
\usepackage{algorithm}
\usepackage{algorithmic}
\usepackage{booktabs}
\usepackage{capt-of}
\usepackage{url}


\DeclareGraphicsExtensions{.pdf,.png,.jpg}

\newcommand{\R}{\mathbb{R}}

\newcommand{\norm}[1]{\left\lVert #1 \right\rVert}

\newcommand{\qules}{\textsc{QULES}}

\theoremstyle{thmstyleone}
\newtheorem{theorem}{Theorem}
\newtheorem{lemma}{Lemma}
\newtheorem{proposition}{Proposition}
\newtheorem{corollary}{Corollary}
\theoremstyle{thmstylethree}
\newtheorem{assumption}{Assumption}

\usepackage{etoolbox}
\makeatletter
\patchcmd{\@maketitle}{\vskip21pt}{\vskip0pt}{}{}
\patchcmd{\@maketitle}{\removelastskip\vskip20pt}{\removelastskip\vskip6pt}{}{}
\patchcmd{\@maketitle}{\removelastskip\vskip24pt}{\removelastskip\vskip12pt}{}{}
\makeatother

\renewcommand{\pacsbullet}{,\ }

\newcommand{\doi}[1]{{\urlstyle{same}\url{https://doi.org/#1}}}

%% file: sections/introduction.tex
\section{Introduction}
\label{sec:introduction}

Derivative-free optimization concerns the minimization of
objective functions for which derivative information is
unavailable or impractical to obtain~\cite{ConnScheinbergVicente2009,AudetHare2017}. Such problems arise
naturally in simulation-driven settings, where each evaluation
may require the execution of a computationally expensive model
and the effective use of function evaluations becomes a central
algorithmic challenge. Derivative-free methods exploit function
values through mechanisms that include constructing local models
from previously evaluated points and generating trial
points using prescribed, adaptive, or randomized
procedures~\cite{LarsonMenickellyWild2019,RiosSahinidis2013,KoldaLewisTorczon2003}.

Convergence analyses of direct-search methods are closely tied to
the geometry of the polling directions. For continuously differentiable objectives,
classical deterministic analyses establish first-order stationarity
results using positive-spanning polling directions that provide sufficient directional
coverage~\cite{Torczon1997}. Randomized direct-search methods replace such deterministic
coverage requirements with probabilistic descent conditions~\cite{GrattonRoyerVicenteZhang2015}.
For locally Lipschitz objectives, the MADS framework instead uses
refining polling directions that become asymptotically dense to obtain
nonsmooth stationarity results~\cite{AudetDennis2006}.

Model-based derivative-free methods construct local linear, quadratic,
or more general surrogate models from sampled function values and
generate trial steps by approximately minimizing these models within
a trust region~\cite{ConnScheinbergVicente2009Global,Powell2002,Powell2004}.
In model-based trust-region frameworks, convergence analyses typically
rely on fully linear or fully quadratic model-accuracy conditions.
These conditions depend on the scale and geometry of the interpolation or
regression set, motivating geometry-improving or self-correcting mechanisms
to maintain suitable poisedness~\cite{ConnScheinbergVicente2008a,ConnScheinbergVicente2008b,ScheinbergToint2010}.

We propose an envelope-search approach for deterministic derivative-free
optimization. Rather than generating a trial point from a polling set or
by minimizing a surrogate model, the proposed Quadratic Upper--Lower Envelope
Search (\qules{}) method determines a candidate through two quadratic functions
with a common center and offset, using their common center as the trial point.
These functions are constructed to bound the sampled objective values from below
and above, while minimizing the curvature gap subject to a trust-region constraint.
After evaluating the trial point, \qules{} updates the incumbent, trust-region
radius, and maintained sample set according to the sufficient-decrease test.

The scalar envelopes use a single curvature parameter for each envelope,
while the diagonal formulation allows the curvature to vary across coordinate
directions. We introduce \qules{}-R to extend this flexibility by constructing
the envelopes in an adaptive orthogonal frame. \qules{}-R estimates curvature
from differences in sampled function values and uses it to update the frame
through pairwise rotations. These updates allow the envelope directions to change
between iterations rather than remain aligned with the coordinate axes.
This mechanism uses points in the sample set and requires no additional
objective evaluations.

The convergence analysis of \qules{} differs from conventional direct-search and
model-based trust-region analyses. The analysis requires neither fully linear
or fully quadratic models nor an optimal solution of the nonconvex envelope
subproblem. Requiring a minimizer would make the convergence guarantees depend
on an outcome the numerical solver may fail to deliver, since a solver can stop
at a feasible point without establishing even local optimality. The proofs rely
instead on feasible quadratic envelopes, which are always available, together
with a bound on the upper-envelope curvature and a sample-geometry condition.
Under these conditions, a hybrid Lyapunov function establishes that the
trust-region radii converge to zero, while an envelope rigidity estimate
provides the basis for proving first-order stationarity.

For deterministic direct-search methods, worst-case first-order complexity bounds of order
$\mathcal{O}(n^2\varepsilon^{-2})$ function evaluations have been established under
standard positive-spanning-set constructions~\cite{Vicente2013}.
Randomized methods like probabilistic descent~\cite{GrattonRoyerVicenteZhang2015}
and Gaussian smoothing~\cite{NesterovSpokoiny2017} reduce the dependence to
$\mathcal O(n\varepsilon^{-2})$ in expectation or with high probability.
Derivative-free trust-region methods using fully linear models achieve a
worst-case evaluation complexity bound of order $\mathcal{O}(n^2\varepsilon^{-2})$
under a suitable scaling of the model-accuracy
constants~\cite{GarmanjaniJudiceVicente2016}. In a geometry-correcting
interpolation framework, standard interpolation-error estimates give an
$\mathcal{O}(n^3\varepsilon^{-2})$ bound. Tighter control of the
interpolation geometry can reduce this bound to
$\mathcal{O}(n^2\varepsilon^{-2})$~\cite{ChaudhryScheinberg2025}.

We establish an $\mathcal{O}(n^2\varepsilon^{-2})$ evaluation complexity 
bound for \qules{} for objectives with Lipschitz-continuous gradients.
This bound holds when the sample set size is $\mathcal{O}(n)$, the geometry
threshold is bounded away from zero, and the remaining bound factors are
uniformly controlled across dimensions. For \qules{}-R, which retains a
sample set of size $\mathcal{O}(n^2)$, the evaluation complexity bound
becomes $\mathcal{O}(n^3\varepsilon^{-2})$.

\noindent\textbf{Contributions.}
The main contributions of this paper are as follows:
\begin{itemize}
  \item We introduce \qules{}, a novel derivative-free method that constructs
  lower and upper quadratic envelopes around the sampled objective values and
  selects their common center as the trial point.

  \item We also introduce \qules{}-R, which adapts the orthogonal envelope frame
  from previously evaluated function values without additional objective
  evaluations. The analysis applies to any sequence of orthogonal frames.

  \item Without any smoothness assumption on $f$, we show that a hybrid
  Lyapunov function decreases at every iteration, so the trust-region radii
  tend to zero when the objective values along the iterates are bounded below.

  \item For continuously differentiable objectives with a compact initial
  lower level set, assuming a uniform upper
  bound on the upper-envelope curvature only when the trust-region radius
  decreases, we prove that
  $\|\nabla f(x_k)\|\to0$ along the whole sequence of iterates, using an
  envelope rigidity estimate and a weighted sample-geometry
  condition.

  \item Additionally, for Lipschitz-continuous gradients, we bound the number
  of evaluations, including geometry repair, by $\mathcal{O}(n^2\varepsilon^{-2})$
  for \qules{} and by $\mathcal{O}(n^3\varepsilon^{-2})$ for \qules{}-R with
  $\mathcal{O}(n^2)$ samples.
\end{itemize}

%% file: sections/problem-setting.tex
\section{Problem setting and notation}
\label{sec:problem-setting}

We consider the unconstrained optimization problem
\begin{equation}
\min_{x\in\mathbb{R}^n} f(x),
\label{eq:problem}
\end{equation}
where $f:\mathbb{R}^n\to\mathbb{R}$ is expensive to evaluate and is 
accessed through a deterministic zeroth-order oracle. 
Given a query point $x\in\mathbb{R}^n$, the oracle returns 
the exact function value $f(x)$. No first- or higher-order derivative 
information is available. We assume that function evaluations dominate
computational cost and count oracle calls. Smoothness assumptions on $f$ are introduced
with the convergence and complexity analyses in
Sect.~\ref{sec:theory}.

At iteration $k$, the algorithm maintains a finite, nonempty sample set\linebreak
$X_k=\{x_{k,1},\ldots,x_{k,m_k}\}$ of evaluated points and their corresponding
objective values. Here $m_k:=|X_k|$, and $\bar m$ denotes a uniform upper
bound on the sample set size, including points added during geometry improvement.
The algorithm chooses the current iterate $x_k$ as an incumbent point satisfying
\begin{equation}
x_k\in\arg\min_{x\in X_k} f(x).
\label{eq:incumbent}
\end{equation}

The envelope subproblem generates a trial point $x_k^+$. We denote 
the displacement from $x_k^+$ to $x_k$ by
\[
p_k:=x_k-x_k^+,
\qquad
\|p_k\|\leq\Delta_k,
\]
where $\Delta_k>0$ is the trust-region radius.

Throughout this paper, $\|\cdot\|$ denotes the Euclidean norm
and $\mathbb R_+$ denotes the set of nonnegative real numbers. For a
continuously differentiable objective, a point $x_\star$ is first-order
stationary if
\[
    \nabla f(x_\star)=0.
\]
Given $\varepsilon>0$, an iterate $x_k$ is
$\varepsilon$-first-order stationary if
\[
    \|\nabla f(x_k)\|\leq\varepsilon.
\]

%% file: sections/envelope-search-algorithm.tex
\section{Envelope Search Algorithm}
\label{sec:algorithm}

\begin{figure}[t]
\centering
\includegraphics[width=0.52\linewidth]{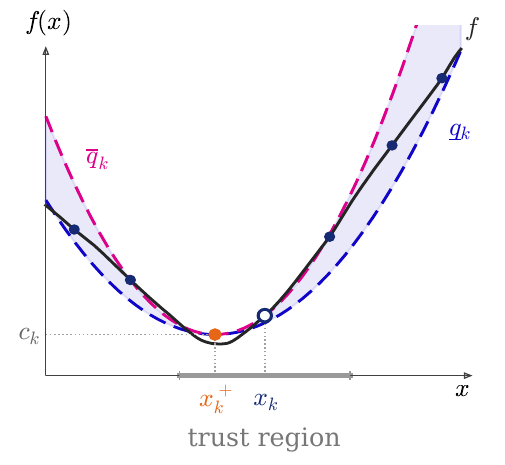}
\caption{Quadratic envelopes with common center $x_k^+$ and a trust region
centered at $x_k$.}
\label{fig:envelope-trial-geometry}
\end{figure}

We define the envelope subproblem, an inner optimization problem that
constructs a pair of quadratic envelopes providing lower and upper bounds on
the sampled objective values.
The envelopes have a common offset $c_k$ and curvature parameters $\alpha_k$
and $\beta_k$. Their common center $x_k^+$ determines the trial point, while
the trust region is centered at the incumbent $x_k$. The algorithm evaluates
the trial point and uses the sufficient-decrease test to update the incumbent
and radius. Figure~\ref{fig:envelope-trial-geometry} illustrates
the geometric construction of the quadratic envelopes and the trial point.

The inner problem is
\begin{equation}
\label{eq:inner-envelope-problem}
\begin{aligned}
\min_{x_k^+,\, c_k,\, \alpha_k,\, \beta_k} \quad & \beta_k - \alpha_k \\
\text{s.t.} \quad 
& y_i \geq \underline{q}_k(x_i), \quad \forall i, \\
& y_i \leq \overline{q}_k(x_i), \quad \forall i, \\
& f(x_k) = c_k + \frac{\beta_k}{2} \|x_k - x_k^+\|^2, \\
& \|x_k - x_k^+\| \leq \Delta_k, \\
& \beta_k \geq \alpha_k \geq 0, \quad \beta_k \geq \beta_{\min}
\end{aligned}
\end{equation}
where the lower and upper quadratic models are defined as
\begin{equation*}
\label{eq:quad_bounds}
\begin{aligned}
\underline{q}_k(x) &= c_k + \frac{\alpha_k}{2} \|x - x_k^+\|^2, \\
\overline{q}_k(x) &= c_k + \frac{\beta_k}{2} \|x - x_k^+\|^2.
\end{aligned}
\end{equation*}

The data points satisfy $x_i\in X_k$ and $y_i=f(x_i)$. The current iterate $x_k$ is maintained as the 
incumbent point in the sample set, so that $f(x_k)\leq f(x_i)$ for all $i$. The quadratic constraints 
$y_i \geq \underline{q}_k(x_i)$ and $y_i \leq \overline{q}_k(x_i)$
ensure that all observed function values lie between the lower and upper
models, so the pair $(\underline{q}_k, \overline{q}_k)$ forms a valid envelope for the data. 
The trust-region constraint $\norm{x_k-x_k^+}\leq\Delta_k$ restricts the
trial center $x_k^+$ to remain close to the current iterate.
The parameter $\beta_{\min}>0$ is a fixed lower bound on the upper-envelope
curvature.

The touching constraint ensures that the upper envelope and the objective
function have the same value at the incumbent $x_k$. Since the upper envelope also bounds the sampled
values from above, we have
\[
  \norm{x_i-x_k^+}^2-\norm{x_k-x_k^+}^2
  \geq \frac{2}{\beta_k}\bigl(f(x_i)-f(x_k)\bigr)
  \geq 0.
\]
Thus the trial point $x_k^+$ lies in the Voronoi cell of $x_k$.

The objective $\beta_k-\alpha_k$ minimizes the gap between the upper and
lower curvature parameters among feasible envelope pairs.
We first show that the inner problem is always feasible.

\begin{proposition} 
\label{prop:inner-feasibility}
Since $x_k$ minimizes $f$ over $X_k$, the inner
problem~\eqref{eq:inner-envelope-problem} is feasible.
\end{proposition}

\begin{proof}
Set $x_k^+=x_k$ and $c_k=f(x_k)$. This choice satisfies the touching and
trust-region constraints, and both envelope constraints hold with equality at $x_k$.
If $X_k=\{x_k\}$, the choice $\alpha_k=0$ and $\beta_k=\beta_{\min}$ is feasible.
Otherwise, at each sample $x_i\ne x_k$, the envelope constraints reduce to
\[
  \alpha_k\leq d_{k,i}\leq\beta_k,
  \qquad\text{where}\qquad
  d_{k,i}:=\frac{2\bigl(f(x_i)-f(x_k)\bigr)}{\norm{x_i-x_k}^2}.
\]
Since $x_k$ is a minimizer of $f$ over $X_k$, we have $d_{k,i}\geq0$. We
choose the largest lower curvature and the smallest upper curvature allowed
by these inequalities and by the curvature floor,
\[
  \alpha_k^{(0)}:=\min_i d_{k,i},
  \qquad
  \beta_k^{(0)}:=\max\Bigl\{\beta_{\min},\,\max_i d_{k,i}\Bigr\},
\]
where the minimum and the maximum are taken over
$x_i\in X_k\setminus\{x_k\}$. Then
$0\leq\alpha_k^{(0)}\leq d_{k,i}\leq\beta_k^{(0)}$ for every such $i$ and
$\beta_k^{(0)}\geq\beta_{\min}$. Hence
$(x_k^+,c_k,\alpha_k,\beta_k)=(x_k,f(x_k),\alpha_k^{(0)},\beta_k^{(0)})$
satisfies all constraints of \eqref{eq:inner-envelope-problem}.
\end{proof}

\suppressfloats[t]
\begin{figure}[t]
\centering
\includegraphics[width=0.88\linewidth]{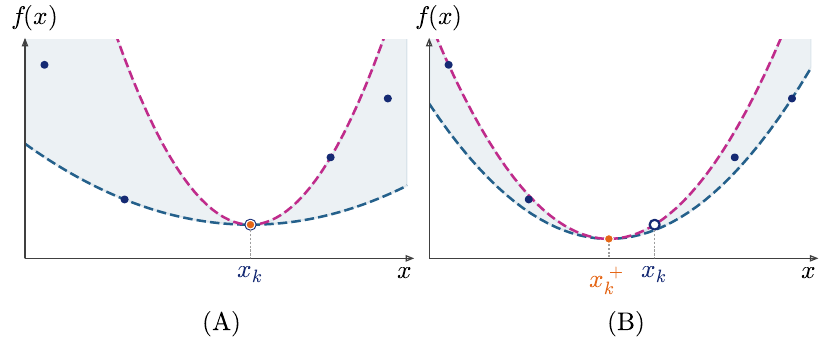}
\caption{Two feasible quadratic envelope pairs for the same sampled values.
Panel (A) has its common center at the incumbent, while panel (B) has
its common center shifted away from the incumbent and a smaller curvature gap.}
\label{fig:envelope-construction}
\end{figure}

Since the inner problem is a nonconvex nonlinear program, obtaining
a numerically reliable solution may be challenging in practice.
We therefore do not require a globally or locally optimal solution.
The nonlinear programming solver is started from the feasible point
constructed in the proof of Proposition~\ref{prop:inner-feasibility}, which
has its common center at the incumbent, as in panel (A) of
Fig.~\ref{fig:envelope-construction}, and can reduce the curvature gap by
moving the common center, as in panel (B). The returned point
$(x_k^+,c_k,\alpha_k,\beta_k)$ is required to be feasible for
\eqref{eq:inner-envelope-problem} and to satisfy
\[
  \beta_k-\alpha_k\leq\beta_k^{(0)}-\alpha_k^{(0)},
\]
which is always possible by Proposition~\ref{prop:inner-feasibility}.
After solving the inner problem~\eqref{eq:inner-envelope-problem},
we evaluate $f(x_k^+)$ and accept the trial point if
\begin{equation}
  f(x_k^+)\leq
  f(x_k)-\eta\frac{\beta_k}{2}\norm{x_k-x_k^+}^2,
  \label{eq:algorithm-acceptance-test}
\end{equation}
where $\eta\in(0,1]$ is the sufficient-decrease parameter.
The upper envelope predicts a decrease of
$\frac{\beta_k}{2}\norm{x_k-x_k^+}^2$. The acceptance test requires
the actual decrease to reach at least a fraction $\eta$ of that amount.

The following proposition shows how far the trial point must move to
achieve a given reduction in the curvature gap. If the step is accepted,
it also gives a lower bound on the decrease in the objective value.

\begin{proposition}
\label{prop:gap-reduction}
Suppose $X_k$ contains a sample point distinct from $x_k$, and let
\[
  r_k:=\min_{x_i\in X_k\setminus\{x_k\}}\norm{x_i-x_k}
\]
be the distance from the incumbent to the nearest other sample point.
For any feasible solution of \eqref{eq:inner-envelope-problem}, let
$g_k:=\beta_k-\alpha_k$ be the curvature gap and let
$G_k:=(\beta_k^{(0)}-\alpha_k^{(0)})-g_k$ be its reduction from the initial
value. If $G_k\geq0$, then
\[
  \norm{x_k^+-x_k}\geq r_k\left(\sqrt{4+\frac{G_k}{\beta_k}}-2\right).
\]
If the trial point also satisfies \eqref{eq:algorithm-acceptance-test}, then
\[
  f(x_k)-f(x_k^+)\geq
  \frac{\eta\beta_k r_k^2}{2}\left(\sqrt{4+\frac{G_k}{\beta_k}}-2\right)^2.
\]
\end{proposition}

\begin{proof}
Set $u:=\norm{x_k^+-x_k}/r_k$. Eliminating $c_k$ by the touching constraint,
expanding the envelope inequalities, and applying the Cauchy--Schwarz inequality
gives, for each $x_i\neq x_k$,
\[
  \alpha_k(1-2u)-g_k u^2\leq d_{k,i}\leq\beta_k(1+2u).
\]
By the definitions of $\alpha_k^{(0)}$ and $\beta_k^{(0)}$ and the floor
$\beta_k\geq\beta_{\min}$, it follows that
\[
  G_k\leq 2(\alpha_k+\beta_k)u+g_k u^2
  \leq\beta_k(4u+u^2).
\]
Solving for $u$ gives the step-length bound. Substitution into
\eqref{eq:algorithm-acceptance-test} gives the objective decrease bound.
\end{proof}

At every iteration, the algorithm adds the evaluated trial point $x_k^+$ to the
maintained dataset before selecting the next iterate. If the trial point satisfies the
sufficient-decrease test \eqref{eq:algorithm-acceptance-test}, it becomes the
reference point for the next geometry-improvement step:
\[
\hat{x}_{k+1}=x_k^+.
\]
Fix a step threshold $\kappa_\Delta\in(0,1)$ and factors $\gamma>1$ and
$\theta_1,\theta_2\in(0,1)$. Each iteration produces one of three outcomes.
A \emph{large-step success} satisfies \eqref{eq:algorithm-acceptance-test}
with $\norm{x_k-x_k^+}\geq\kappa_\Delta\Delta_k$. The step is large relative
to the trust region, so the radius is enlarged to
$\Delta_{k+1}=\gamma\Delta_k$ to explore a larger neighborhood.
A \emph{small-step success} satisfies \eqref{eq:algorithm-acceptance-test}
with $\norm{x_k-x_k^+}<\kappa_\Delta\Delta_k$. The step is small relative to
the trust region, so the radius is reduced to
$\Delta_{k+1}=\theta_1\Delta_k$ to refine the local search.
A \emph{failed trial} does not satisfy \eqref{eq:algorithm-acceptance-test},
and the radius is reduced to $\Delta_{k+1}=\theta_2\Delta_k$.

A trial point that fails the sufficient-decrease test may still have a lower
objective value than the current iterate, because the decrease may not be
large enough to satisfy the test. In this case, an incumbent of the updated
sample set is selected as the reference point for the next geometry-improvement
step,
\[
  \hat{x}_{k+1}\in\arg\min_{x\in X_{k+1}} f(x).
\]
Thus, $x_k^+$ may still become the reference point. The subsequent geometry
improvement preserves $f(x_{k+1})\leq f(x_k)$.

Fix a pruning factor $\rho_s\geq1$. After an accepted trial, the algorithm
prunes points farther than $\rho_s\Delta_{k+1}$ from the new incumbent without
reducing the sample set size below a prescribed minimum. After an unsuccessful
trial that adds a new point, it removes the farthest old point if the sample set
exceeds this minimum. In both cases, it enforces a fixed capacity by removing
the farthest points.

Before each inner solve, a geometry-improvement step checks the sample
geometry and adds points if necessary. Let $\hat{x}_k$ denote the incumbent
at the start of this step, which is used as a fixed reference point. Define
\[
  w_{k,i}:=\frac{\Delta_k^2}{\Delta_k^2+\|x_{k,i}-\hat{x}_k\|^2},
  \qquad
  Z_k:=\left[\,\frac{w_{k,i}}{\Delta_k}(x_{k,i}-\hat{x}_k)\,\right]_{i=1}^{m_k}.
\]
The sample set is prepared to satisfy the weighted geometry condition
\begin{equation}
  \min_{\|v\|=1}\|Z_k^\top v\|\geq\sigma_s,
  \label{eq:sample-geometry}
\end{equation}
where $\sigma_s$ is the geometry threshold. If the condition is not satisfied,
the algorithm repairs the set about the fixed reference point $\hat{x}_k$: it
selects a least-represented direction $u$ from the singular value decomposition
of $Z_k$, evaluates $f(\hat{x}_k+\Delta_k u)$, and adds the point to $X_k$.
It repeats the test with the updated matrix until the condition holds, keeping
$\hat{x}_k$ and $\Delta_k$ fixed and retaining all samples during repair.
After geometry improvement, it selects the incumbent $x_k$ for the inner solve.
If no strict improvement is found, it sets $x_k=\hat{x}_k$.
The weights and $Z_k$ continue to use $\hat{x}_k$ as their reference point.

At each iteration, the algorithm combines the trust-region update mechanism
with the adaptive sample set management described above. The resulting
procedure is summarized in Algorithm~\ref{alg:qules}.

\begin{algorithm}[t]
\caption{\qules{} iteration}
\label{alg:qules}
\begin{algorithmic}[1]
\STATE Choose $\Delta_0>0$, $\beta_{\min}>0$, $\eta\in(0,1]$, $\theta_1,\theta_2\in(0,1)$,
$\gamma>1$, and $\kappa_\Delta\in(0,1)$.
\STATE Set pruning factor $\rho_s\geq1$ and geometry threshold $\sigma_s\in(0,1/2)$.
\STATE Initialize a finite evaluated data set $X_0$.
\STATE Choose $\hat{x}_0\in\arg\min_{x\in X_0} f(x)$.
\FOR{$k=0,1,2,\ldots$}
  \IF{the geometry condition fails at $\hat{x}_k$}
    \STATE Repair geometry about $\hat{x}_k$ at radius $\Delta_k$.
  \ENDIF
  \STATE Choose $x_k\in\arg\min_{x\in X_k} f(x)$, keeping $\hat{x}_k$ on ties.
  \STATE Find $(x_k^+,c_k,\alpha_k,\beta_k)$ feasible for \eqref{eq:inner-envelope-problem} with gap at most the starting gap.
  \STATE Evaluate $f(x_k^+)$.
  \STATE Set $\widehat X_{k+1}=X_k\cup\{x_k^+\}$.
  \IF{$f(x_k^+)\leq f(x_k)-\eta \frac{\beta_k}{2}
  \norm{x_k-x_k^+}^2$}
    \IF{$\norm{x_k-x_k^+}\geq\kappa_\Delta\Delta_k$}
      \STATE Set $\Delta_{k+1}=\gamma\Delta_k$.
    \ELSE
      \STATE Set $\Delta_{k+1}=\theta_1\Delta_k$.
    \ENDIF
    \STATE Set $X_{k+1}=\widehat X_{k+1}$.
    \STATE Set $\hat{x}_{k+1}=x_k^+$.
    \STATE Prune $X_{k+1}$ about $\hat{x}_{k+1}$ at radius $\rho_s\Delta_{k+1}$.
    \STATE Enforce capacity using reference $\hat{x}_{k+1}$.
  \ELSE
    \STATE Set $\Delta_{k+1}=\theta_2\Delta_k$.
    \STATE Set $X_{k+1}=\widehat X_{k+1}$.
    \STATE Remove the farthest old point from $x_k$, subject to the minimum size.
    \STATE Enforce capacity using reference $x_k$.
    \STATE Set $\hat{x}_{k+1}\in\arg\min_{x\in X_{k+1}} f(x)$.
  \ENDIF
\ENDFOR
\end{algorithmic}
\end{algorithm}

%% file: sections/convergence-analysis.tex
\section{Theoretical analysis}
\label{sec:theory}

At each iteration, Algorithm~\ref{alg:qules} ensures that the sample set
satisfies the geometry condition~\eqref{eq:sample-geometry} and contains at
least one point distinct from the incumbent. As shown in
Proposition~\ref{prop:inner-feasibility}, the inner
problem~\eqref{eq:inner-envelope-problem} always has a feasible solution.
We first show that the trust-region radii converge to zero, then establish
first-order stationarity and an evaluation complexity bound. The latter
two results use the following curvature bound.

\begin{assumption}
\label{ass:bounded-curvature}
There exists a finite constant $\beta_{\max}>0$ such that
$\beta_k\leq\beta_{\max}$ whenever $\Delta_{k+1}<\Delta_k$.
\end{assumption}

Uniform boundedness of the model Hessian at every iteration is a standard
assumption in derivative-free trust-region analysis
\cite[Assumption~5.2]{ConnScheinbergVicente2009Global}.
For the upper envelope, whose Hessian is $\beta_k I$ with
$\norm{\beta_k I}_2=\beta_k$, Assumption~\ref{ass:bounded-curvature} requires
this bound only on iterations in which the trust-region radius decreases.
The curvature parameter $\beta_k$ is a decision variable in
\eqref{eq:inner-envelope-problem}, so its value is determined by the inner
solver, and the algorithm has no mechanism for enforcing an upper bound on it.
Unlike the lower bound, the upper bound cannot be imposed as a constraint of
\eqref{eq:inner-envelope-problem}, because the inner problem could then be
infeasible. Proposition~\ref{prop:bounded-upper-curvature} provides sufficient
conditions under which this bound holds.

\begin{proposition}
\label{prop:bounded-upper-curvature}
Let $(x_k^+,c_k,\alpha_k,\beta_k)$ be feasible for
\eqref{eq:inner-envelope-problem} and satisfy
$\beta_k-\alpha_k\leq\beta_k^{(0)}-\alpha_k^{(0)}$, where
$\alpha_k^{(0)}$ and $\beta_k^{(0)}$ are the initial curvatures from
Proposition~\ref{prop:inner-feasibility}. Define
\[
  \lambda_k:=
  \frac{\norm{x_k^+-x_k}}{
    \displaystyle\max_{x_i\in X_k\setminus\{x_k\}}\norm{x_k^+-x_i}}.
\]
If $\lambda_k<1$, then
\[
  \beta_k\leq\frac{3\beta_k^{(0)}}{1-\lambda_k}.
\]
Consequently, if there are constants $B>0$ and $\bar\lambda\in[0,1)$
such that $\beta_k^{(0)}\leq B$ and $\lambda_k\leq\bar\lambda$ whenever
$\Delta_{k+1}<\Delta_k$, then Assumption~\ref{ass:bounded-curvature} holds
with $\beta_{\max}=3B/(1-\bar\lambda)$.
\end{proposition}

\begin{proof}
Choose a sample $x_i\neq x_k$ attaining the maximum in the denominator.
Set $t=\norm{x_k^+-x_i}$, $\delta=\norm{x_k^+-x_k}$, and
$g_k=\beta_k-\alpha_k$. Since $x_k^+$ lies in the Voronoi cell of $x_k$,
we have $t\geq\delta$, and $t>0$ since $x_i\neq x_k$. Eliminating $c_k$ via the touching
constraint and using the lower-envelope constraint at $x_i$, we obtain
\[
  \beta_k(t^2-\delta^2)\leq 2\bigl(f(x_i)-f(x_k)\bigr)+g_k t^2.
\]
By the definition of $\beta_k^{(0)}$ and the triangle inequality,
\[
  2\bigl(f(x_i)-f(x_k)\bigr)
  \leq\beta_k^{(0)}\norm{x_i-x_k}^2
  \leq\beta_k^{(0)}(1+\lambda_k)^2t^2.
\]
Since $g_k\leq\beta_k^{(0)}$ and $\delta=\lambda_k t$, dividing by
$t^2(1-\lambda_k^2)$ gives
\[
  \beta_k\leq\beta_k^{(0)}
  \frac{(1+\lambda_k)^2+1}{1-\lambda_k^2}
  \leq\frac{3\beta_k^{(0)}}{1-\lambda_k}.
  \qedhere
\]
\end{proof}

On iterations that reduce the trust-region radius,
Proposition~\ref{prop:bounded-upper-curvature} gives a uniform bound on
$\beta_k$ provided that the initial curvatures remain bounded and the trial
steps satisfy the stated geometric condition. A zero step gives
$\lambda_k=0$. If the step is at most $q<1/2$ times the distance from the
incumbent to its nearest other sample, then $\lambda_k\leq q/(1-q)$.
Proposition~\ref{prop:bounded-upper-curvature} also covers steps that are
not small relative to the nearest sample. It suffices that some sample lies
at distance at least $c\norm{x_k^+-x_k}$ from the trial point for a fixed
$c>1$, since then $\lambda_k\leq 1/c$. For the initial curvature bound, one may take
$B=\max\{\beta_{\min},2C\}$ if $f(x_i)-f(x_k)\leq C\norm{x_i-x_k}^2$
for all samples on those iterations, with the same $C$ throughout.

\subsection{Hybrid Lyapunov analysis}
\label{subsec:phase1}

At each iteration, the trust-region mechanism described in
Sect.~\ref{sec:algorithm} produces one of three outcomes: a failed trial, a
large-step success, or a small-step success. Since the radius is enlarged only
after a large-step success and reduced otherwise, its asymptotic behavior is
not immediate. This mode-dependent update motivates a hybrid Lyapunov function
defined on the joint state
\[
  z_k:=(x_k,\Delta_k)\in\R^n\times\R_+.
\]
For $\omega>0$, define
\[
V_k:=V(z_k):=f(x_k)+\omega\Delta_k^2.
\]
Using the trial displacement $p_k$, the acceptance condition
\eqref{eq:algorithm-acceptance-test} gives
\begin{equation}
f(x_k)-f(x_k^+)
\geq
\eta\frac{\beta_k}{2}\norm{p_k}^2
\label{eq:successful-decrease-phase1}
\end{equation}
whenever the trial point is accepted.

\begin{proposition} 
\label{prop:lyapunov-decrease}
Let the iterates be updated by Algorithm~\ref{alg:qules}. If
\begin{equation}
  0<\omega<
  \frac{\eta\beta_{\min}\kappa_\Delta^2}{2(\gamma^2-1)},
  \label{eq:omega-choice}
\end{equation}
then there exists $\mu>0$ such that
\begin{equation}
  V_k-V_{k+1}\geq\mu\Delta_k^2
  \qquad\text{for all }k\geq0.
  \label{eq:uniform-lyapunov-decrease}
\end{equation}
\end{proposition}

\begin{proof}
Suppose first that iteration $k$ is a failed trial. Since
$\Delta_{k+1}=\theta_2\Delta_k$, we have
\begin{align*}
V_k-V_{k+1}
&=f(x_k)-f(x_{k+1})
  +\omega(1-\theta_2^2)\Delta_k^2 \\
&\geq \omega(1-\theta_2^2)\Delta_k^2.
\end{align*}
Suppose next that iteration $k$ is a large-step success. Then
$\Delta_{k+1}=\gamma\Delta_k$ and
$\norm{p_k}\geq\kappa_\Delta\Delta_k$. Using
\eqref{eq:successful-decrease-phase1} and
$\beta_k\geq\beta_{\min}$ gives
\begin{align*}
V_k-V_{k+1}
&=f(x_k)-f(x_{k+1})
  -\omega(\gamma^2-1)\Delta_k^2 \\
&\geq f(x_k)-f(x_k^+)
  -\omega(\gamma^2-1)\Delta_k^2 \\
&\geq \eta\frac{\beta_k}{2}\norm{p_k}^2
  -\omega(\gamma^2-1)\Delta_k^2 \\
&\geq
\left(
\eta\frac{\beta_{\min}}{2}\kappa_\Delta^2
-\omega(\gamma^2-1)
\right)\Delta_k^2.
\end{align*}
The coefficient in the final line is positive by
\eqref{eq:omega-choice}.
Finally, suppose that iteration $k$ is a small-step success. Since
$\Delta_{k+1}=\theta_1\Delta_k$, it follows from
\eqref{eq:successful-decrease-phase1} that
\begin{align*}
V_k-V_{k+1}
&=f(x_k)-f(x_{k+1})
  +\omega(1-\theta_1^2)\Delta_k^2 \\
&\geq f(x_k)-f(x_k^+)
  +\omega(1-\theta_1^2)\Delta_k^2 \\
&\geq \omega(1-\theta_1^2)\Delta_k^2.
\end{align*}
Therefore, \eqref{eq:uniform-lyapunov-decrease} holds with
\begin{equation*}
\mu:=
\min\left\{
\omega(1-\theta_1^2),
\omega(1-\theta_2^2),
\eta\frac{\beta_{\min}}{2}\kappa_\Delta^2
-\omega(\gamma^2-1)
\right\}>0.
\end{equation*}
The first two terms are positive since $\theta_1,\theta_2\in(0,1)$,
and the third is positive by \eqref{eq:omega-choice}.
\end{proof}

Proposition~\ref{prop:lyapunov-decrease} shows that the Lyapunov function decreases at 
each iteration. To bound the accumulated decrease, we now require
the objective values along the generated iterates to be bounded below.

\begin{theorem} 
\label{thm:radius-decay}
Suppose $f(x_k) \geq f_{\min} > -\infty$ for all $k\geq0$. For any analysis
constant $\omega$ satisfying~\eqref{eq:omega-choice},
\[
  \sum_{k=0}^{\infty}\Delta_k^2<\infty,
  \qquad\text{and hence}\qquad
  \Delta_k\to0.
\]
\end{theorem}

\begin{proof}
By Proposition~\ref{prop:lyapunov-decrease},
\[
\mu\sum_{k=0}^{N}\Delta_k^2 \le \sum_{k=0}^{N}(V_k-V_{k+1}) = V_0-V_{N+1}
\]
for every $N\geq0$. Since
$V_{N+1}=f(x_{N+1})+\omega\Delta_{N+1}^2\geq f_{\min}$, it follows that
\[
\mu\sum_{k=0}^{N}\Delta_k^2 \le V_0-f_{\min}.
\]
Hence
\[
\sum_{k=0}^{\infty}\Delta_k^2<\infty.
\]
The nonnegativity of the terms implies $\Delta_k^2\to0$, and therefore $\Delta_k\to0$.
\end{proof}

This section establishes the asymptotic stability of the trust-region mechanism in the 
sense that the radii satisfy $\Delta_k \to 0$. Note that this conclusion is obtained 
without any differentiability or smoothness assumption on $f$, relying only on the 
sufficient-decrease mechanism, the curvature floor $\beta_k\geq\beta_{\min}$ in \eqref{eq:inner-envelope-problem}, 
and lower boundedness along the generated iterates. However, this conclusion alone provides 
no information about the first-order stationarity of the iterates.

\subsection{Envelope rigidity}
\label{subsec:phase2}

The first-order analysis links the trust-region mechanism to the local information
provided by the envelope inequalities. This connection allows the behavior of the algorithm
to be related to first-order stationarity. It is sufficient to impose mild regularity 
on the objective.

\begin{assumption} 
\label{ass:objective-regularity}
The function $f:\R^n\to\R$ is continuously differentiable, and the initial
lower level set
$\mathcal L(x_0):=\{x\in\mathbb R^n : f(x)\leq f(x_0)\}$ is compact.
\end{assumption}

The compact level set contains all incumbents, and continuous differentiability gives 
a first-order Taylor remainder that is uniform near this set. To connect this local 
information to the envelope inequalities, we first establish a geometric bound for 
the weighted sample set.

\begin{lemma}
\label{lem:sample-geometry}
Let $X_k$ be the prepared sample set at iteration $k$, where $m_k\leq\bar m$,
and let $z_{k,1},\ldots,z_{k,m_k}$ be the columns of $Z_k$. Then
\[
  \|v\|
  \leq\frac{\sqrt{\bar m}}{\sigma_s}
  \max_{1\leq i\leq m_k}|z_{k,i}^{\top}v|
  \qquad\text{for every }v\in\mathbb R^n.
\]
\end{lemma}

\begin{proof}
For $v\neq0$, applying the geometry condition
\eqref{eq:sample-geometry} to $v/\|v\|$ gives
\[
  \sigma_s\|v\|\leq\|Z_k^\top v\|.
\]
Since the columns of $Z_k$ are $z_{k,1},\ldots,z_{k,m_k}$ and $m_k\leq\bar m$,
\[
  \|Z_k^\top v\|
  =\left(\sum_{i=1}^{m_k}|z_{k,i}^\top v|^2\right)^{1/2}
  \leq\sqrt{\bar m}\max_{1\leq i\leq m_k}|z_{k,i}^\top v|.
\]
Combining these inequalities and dividing by $\sigma_s$ proves the result.
The case $v=0$ is immediate.
\end{proof}

Having established a bound based on the weighted sample geometry, we now use
the envelope inequalities to bound the weighted function differences.

\begin{lemma}
\label{lem:envelope-rigidity}
For every prepared sample set $X_k$,
\[
  w_{k,i}\left|f(x_{k,i})-f(x_k)\right|
  \leq\frac52\beta_k\Delta_k^2
  \qquad\text{for }i=1,\ldots,m_k.
\]
\end{lemma}

\begin{proof}
Fix $k$ and let $s_i:=x_{k,i}-x_k$ and $r_i:=\|x_{k,i}-\hat{x}_k\|$.
Since $x_k$ is an incumbent of $X_k$, we have $f(x_{k,i})-f(x_k)\geq0$.
The upper envelope and touching condition give
\[
\begin{aligned}
  f(x_{k,i})-f(x_k)
  &\leq\beta_k p_k^\top s_i+\frac{\beta_k}{2}\|s_i\|^2\\
  &\leq\beta_k\left(\Delta_k\|s_i\|+\frac12\|s_i\|^2\right).
\end{aligned}
\]
The incumbent $x_k$ is either the fixed reference point $\hat{x}_k$ or a new
point evaluated at distance $\Delta_k$ from it during geometry improvement.
Thus $\|x_k-\hat{x}_k\|\leq\Delta_k$, and the triangle inequality gives
$\|s_i\|\leq r_i+\Delta_k$. Therefore,
\[
\begin{aligned}
  \left|f(x_{k,i})-f(x_k)\right|
  &\leq\beta_k
  \left(\Delta_k(r_i+\Delta_k)+\frac12(r_i+\Delta_k)^2\right)\\
  &\leq\frac52\beta_k(\Delta_k^2+r_i^2),
\end{aligned}
\]
where the last inequality follows from
\[
  \frac52(\Delta_k^2+r_i^2)
  -\Delta_k(r_i+\Delta_k)-\frac12(r_i+\Delta_k)^2
  =(\Delta_k-r_i)^2+r_i^2\geq0.
\]
Multiplying by $w_{k,i}$ and using
$w_{k,i}(\Delta_k^2+r_i^2)=\Delta_k^2$ proves the result.
\end{proof}

Lemma~\ref{lem:envelope-rigidity} shows that the weighted differences between
the sampled objective values and $f(x_k)$ are bounded by
$\frac52\beta_k\Delta_k^2$. This estimate follows directly from the upper
envelope and touching condition.
Assumption~\ref{ass:objective-regularity} is then used to relate these
function-value differences to the first-order variation of $f$ along
the sample displacements.

\begin{lemma}
\label{lem:directional-control}
Suppose Assumptions~\ref{ass:bounded-curvature}
and~\ref{ass:objective-regularity} hold.
If $\Delta_k\to0$, then, along iterations with $\Delta_{k+1}<\Delta_k$,
\[
  \max_{1\leq i\leq m_k}
  |z_{k,i}^\top\nabla f(x_k)|\longrightarrow0.
\]
\end{lemma}

\begin{proof}
Let $K=\mathcal L(x_0)$, $G=\max_{x\in K}\|\nabla f(x)\|$, and let
$\varpi$ be a nondecreasing modulus of continuity of $\nabla f$ on
$K+\overline B(0,1)$. For each sample point $x_{k,i}$, let
$s_i=x_{k,i}-x_k$, $r_i=\|x_{k,i}-\hat{x}_k\|$, and $e_k=\hat{x}_k-x_k$.
If $r_i\leq\sqrt{\Delta_k}$, then, for sufficiently large $k$,
\[
  \|e_k\|\leq\Delta_k,
  \qquad
  \|s_i\|\leq r_i+\|e_k\|\leq2\sqrt{\Delta_k}\leq1,
\]
so the segments $[x_k,x_{k,i}]$ and $[x_k,\hat{x}_k]$ lie in
$K+\overline B(0,1)$.
Since $x_{k,i}-\hat{x}_k=s_i-e_k$, first-order Taylor expansions of
$f(x_{k,i})$ and $f(\hat{x}_k)$ about $x_k$ give
\[
  \nabla f(x_k)^{\top}(s_i-e_k)
  =\bigl(f(x_{k,i})-f(x_k)\bigr)-\bigl(f(\hat{x}_k)-f(x_k)\bigr)
   -R_i+R_e,
\]
where the remainders satisfy
\[
  |R_i|\leq\|s_i\|\,\varpi(2\sqrt{\Delta_k}),
  \qquad
  |R_e|\leq\|e_k\|\,\varpi(2\sqrt{\Delta_k}).
\]
Multiply by $w_{k,i}/\Delta_k$ and take absolute values.
Recall that $w_{k,i}\leq1$.
For the two envelope terms, Lemma~\ref{lem:envelope-rigidity} applied to
$x_{k,i}$, and to $\hat{x}_k$ with weight $1$, gives
\[
  \frac{w_{k,i}}{\Delta_k}\,\bigl|f(x_{k,i})-f(x_k)\bigr|
  \leq\frac52\beta_k\Delta_k,
  \qquad
  \frac{w_{k,i}}{\Delta_k}\,\bigl|f(\hat{x}_k)-f(x_k)\bigr|
  \leq\frac52\beta_k\Delta_k.
\]
For the two remainder terms, $\|s_i\|\leq r_i+\Delta_k$ and
\[
  \frac{w_{k,i}r_i}{\Delta_k}
  =\frac{\Delta_k r_i}{\Delta_k^2+r_i^2}\leq\frac12
\]
give
\[
  \frac{w_{k,i}\|s_i\|}{\Delta_k}\leq\frac32,
  \qquad
  \frac{w_{k,i}\|e_k\|}{\Delta_k}\leq1.
\]
Altogether,
\[
  |z_{k,i}^{\top}\nabla f(x_k)|
  \leq5\beta_k\Delta_k+\frac52\,\varpi(2\sqrt{\Delta_k}).
\]
If $r_i>\sqrt{\Delta_k}$, then
\[
  |z_{k,i}^{\top}\nabla f(x_k)|
  \leq G\frac{\Delta_k r_i}{\Delta_k^2+r_i^2}
  \leq G\sqrt{\Delta_k}.
\]
By Assumption~\ref{ass:bounded-curvature},
\mbox{$\beta_k\Delta_k\leq\beta_{\max}\Delta_k\to0$} on iterations that reduce
the radius, so both bounds tend to zero uniformly over the samples.
\end{proof}

Together with Lemma~\ref{lem:sample-geometry}, the bounds in the proof of
Lemma~\ref{lem:directional-control} show that, once the radius is
sufficiently small, it can be reduced only where the gradient is also small.
Specifically, for every $\varepsilon>0$ there exists $\delta_\varepsilon>0$
such that every iteration with $\|\nabla f(x_k)\|\geq\varepsilon$ and
$\Delta_k<\delta_\varepsilon$ is a large-step success.

\begin{theorem} 
\label{thm:c1-first-order-convergence}
Let the iterates be generated by Algorithm~\ref{alg:qules}.
Suppose Assumptions~\ref{ass:bounded-curvature}
and~\ref{ass:objective-regularity} hold. Then
\[
  \norm{\nabla f(x_k)}\to0.
\]
\end{theorem}

\begin{proof}
All incumbents lie in the compact level set $\mathcal L(x_0)$, on which
$f$ is bounded below. Theorem~\ref{thm:radius-decay} therefore gives
$\Delta_k\to0$. Along iterations that reduce the radius, Lemma~\ref{lem:directional-control} gives
\[
  \max_{1\leq i\leq m_k}
  |z_{k,i}^{\top}\nabla f(x_k)|\longrightarrow0.
\]
Applying Lemma~\ref{lem:sample-geometry} with $v=\nabla f(x_k)$ gives
\[
  \|\nabla f(x_k)\|
  \leq\frac{\sqrt{\bar m}}{\sigma_s}
  \max_{1\leq i\leq m_k}|z_{k,i}^{\top}\nabla f(x_k)|.
\]
Thus $\|\nabla f(x_k)\|\to0$ along iterations that reduce the radius.
For any $k$, let $\ell$ be the first such iteration at or after $k$.
It exists because otherwise the radius would grow by a factor $\gamma>1$
at every subsequent iteration, contradicting $\Delta_k\to0$.
For $k\leq j<\ell$, the accepted trial is the next reference point,
$\hat x_{j+1}=x_j^+$. Geometry improvement moves the incumbent by at most
$\Delta_{j+1}$ from this reference point. Since $\Delta_{j+1}=\gamma\Delta_j$,
\[
  \|x_\ell-x_k\|
  \leq\sum_{j=k}^{\ell-1}(\Delta_j+\Delta_{j+1})
  =\frac{1+\gamma}{\gamma-1}(\Delta_\ell-\Delta_k).
\]
As $k\to\infty$, we have $\ell\geq k$, so $\Delta_\ell\to0$ and
$\|\nabla f(x_\ell)\|\to0$. Uniform continuity of $\nabla f$ on
$\mathcal L(x_0)$ now gives $\|\nabla f(x_k)\|\to0$ along the whole sequence.
\end{proof}

The upper envelope curvature $\beta_k$ need not remain bounded for the
first-order stationarity result to hold. On iterations that reduce the trust-region
radius, it suffices that $\beta_k=o(1/\Delta_k)$, or equivalently, that
$\beta_k\Delta_k\to0$.

\subsection{Evaluation complexity}
\label{subsec:phase3}
The convergence result in Theorem~\ref{thm:c1-first-order-convergence} establishes that 
the generated iterates are asymptotically first-order stationary. However, this result 
does not quantify how many objective evaluations are required before an iterate with
$\|\nabla f(x_k)\| \le \varepsilon$ is obtained. We first bound the cost
of geometry repair and then relate the gradient norm to the trust-region
radius.

\begin{lemma}
\label{lem:geometry-repair}
For $0<\sigma_s<1/2$, the geometry improvement procedure satisfies the weighted
geometry condition \eqref{eq:sample-geometry} after at most $n$ additional
objective evaluations.
\end{lemma}

\begin{proof}
Let $Z$ be the current weighted sample matrix and $A=ZZ^\top$. If the geometry
condition fails, the procedure selects a unit eigenvector $u$ of $A$ with
\[
  \lambda_{\min}(A)=\|Z^\top u\|^2<\sigma_s^2.
\]
The new point $y=\hat{x}_k+\Delta_k u$ has weight $1/2$ and column $u/2$.
It is not already in the sample set, since otherwise
$\|Z^\top u\|\geq1/2>\sigma_s$. Appending the new column gives
\[
  Z_{\mathrm{new}}Z_{\mathrm{new}}^\top=A+\frac14uu^\top.
\]
Since $u$ is an eigenvector of $A$, this rank-one update raises the eigenvalue
$\lambda_{\min}(A)$ by $1/4$, to a value above $\sigma_s^2$, and leaves the
eigenvalues on the orthogonal complement of $u$ unchanged. Each additional
evaluation therefore reduces the number of eigenvalues below $\sigma_s^2$ by
one. Since $A$ has $n$ eigenvalues, the condition \eqref{eq:sample-geometry}
holds after at most $n$ evaluations.
\end{proof}

Here, we strengthen the regularity assumptions 
by requiring $\nabla f$ to be Lipschitz continuous, which provides a uniform quadratic
bound on the first-order Taylor remainder.

\begin{lemma}
\label{lem:gradient-radius-bound}
Suppose $\nabla f$ is $L$-Lipschitz continuous on a set containing the segments
from each incumbent to its sample points. Then
\[
  \norm{\nabla f(x_k)}\leq\frac{\sqrt{\bar m}}{\sigma_s}
  \left(5\beta_k+\frac32L\right)\Delta_k
  \qquad\text{for every }k\geq0.
\]
\end{lemma}

\begin{proof}
Fix $k$ and let
$s_i=x_{k,i}-x_k$, $e_k=\hat{x}_k-x_k$, and
$r_i=\|x_{k,i}-\hat{x}_k\|$.
Taylor expansion along $[x_k,x_{k,i}]$ and $[x_k,\hat{x}_k]$ gives
\[
\begin{aligned}
  |z_{k,i}^\top\nabla f(x_k)|
  &\leq\frac{w_{k,i}}{\Delta_k}
  \left(|f(x_{k,i})-f(x_k)|+|f(\hat{x}_k)-f(x_k)|\right)\\
  &\quad+\frac{Lw_{k,i}}{2\Delta_k}
  \left(\|s_i\|^2+\|e_k\|^2\right).
\end{aligned}
\]
Lemma~\ref{lem:envelope-rigidity}, applied also to $\hat{x}_k$ whose weight
is $1$, bounds the first term by $5\beta_k\Delta_k$.
Moreover, $\|e_k\|\leq\Delta_k$ and $\|s_i\|\leq r_i+\Delta_k$, so
\[
  \frac{w_{k,i}}{\Delta_k}\|s_i\|^2\leq2\Delta_k,
  \qquad
  \frac{w_{k,i}}{\Delta_k}\|e_k\|^2\leq\Delta_k.
\]
Here we use $(r_i+\Delta_k)^2\leq2(r_i^2+\Delta_k^2)$ and
$w_{k,i}(r_i^2+\Delta_k^2)=\Delta_k^2$. Combining these bounds gives
\[
  |z_{k,i}^\top\nabla f(x_k)|
  \leq\left(5\beta_k+\frac32L\right)\Delta_k.
\]
Applying Lemma~\ref{lem:sample-geometry} with $v=\nabla f(x_k)$ gives
\[
  \|\nabla f(x_k)\|
  \leq\frac{\sqrt{\bar m}}{\sigma_s}\max_i|z_{k,i}^\top\nabla f(x_k)|
  \leq\frac{\sqrt{\bar m}}{\sigma_s}\left(5\beta_k+\frac32L\right)\Delta_k,
\]
which proves the result.
\end{proof}

We can now establish the corresponding function-evaluation complexity bound.
The bound includes all objective evaluations used for initialization, one
trial evaluation per iteration, and at most $n$ additional evaluations per
iteration for geometry improvement.

\begin{theorem} 
\label{thm:function-evaluation-complexity}
Suppose Assumption~\ref{ass:bounded-curvature} holds,
$\nabla f$ is $L$-Lipschitz continuous on the sample segments, and
$f(x_k)\geq f_{\min}>-\infty$ for every $k\geq0$.
Let $\omega$ satisfy~\eqref{eq:omega-choice}, let $\mu$ be the constant in
Proposition~\ref{prop:lyapunov-decrease}, and define
\[
  M:=\frac{\sqrt{\bar m}}{\sigma_s}\left(5\beta_{\max}+\frac32L\right).
\]
Then, for every $\varepsilon>0$, Algorithm~\ref{alg:qules} obtains an iterate
satisfying
\[
  \norm{\nabla f(x_k)}\leq\varepsilon
\]
within
\[
  N_f(\varepsilon)
  \leq N_f^{(0)}
  +(n+1)\frac{V_0-f_{\min}}{\mu}
  \max\!\left\{\frac{1}{\Delta_0^2},
  \frac{M^2}{\min\{\theta_1,\theta_2\}^2\varepsilon^2}\right\}
\]
objective evaluations, where $N_f^{(0)}$ counts all evaluations used to
prepare the initial sample set and select $x_0$.
\end{theorem}

\begin{proof}
Fix $\varepsilon>0$. If $\|\nabla f(x_k)\|>\varepsilon$ and the radius
is reduced, Lemma~\ref{lem:gradient-radius-bound} and
Assumption~\ref{ass:bounded-curvature} give $\varepsilon<M\Delta_k$.
The updated radius therefore exceeds $\min\{\theta_1,\theta_2\}\varepsilon/M$.
On every other iteration the radius increases. Hence, by induction, if
none of the first $N$ iterates meets the target, then
\[
  \Delta_k\geq\min\!\left\{\Delta_0,\frac{\min\{\theta_1,\theta_2\}\varepsilon}{M}\right\},
  \qquad 0\leq k<N.
\]
Summing Proposition~\ref{prop:lyapunov-decrease} over these iterations gives
\[
  N\mu\min\!\left\{\Delta_0^2,
  \frac{\min\{\theta_1,\theta_2\}^2\varepsilon^2}{M^2}\right\}
  \leq V_0-V_N\leq V_0-f_{\min}.
\]
This bounds the number of iterations before the first target iterate.
Each iteration costs at most $n+1$ evaluations, including geometry
improvement. Adding the initialization cost gives the stated bound.
\end{proof}

To describe the dimension dependence of the evaluation bound, we take the
Lipschitz constant $L$, the curvature bound $\beta_{\max}$, the geometry
threshold $\sigma_s$, the initial Lyapunov gap, and the remaining algorithmic
parameters to be independent of $n$. The constant $M$ in
Theorem~\ref{thm:function-evaluation-complexity} then satisfies
$M^2=\mathcal O(\bar m)$, and each iteration costs at most $n+1$ evaluations, so
\[
  N_f(\varepsilon)
  \leq N_f^{(0)}+\mathcal O(n\bar m\varepsilon^{-2}),
  \qquad 0<\varepsilon\leq1.
\]
The factor $\bar m$ comes from the sample geometry bound and the factor $n$
from the evaluation cost per iteration. When $\bar m$ scales linearly with $n$
and the initial sample set is constructed using $\mathcal O(n)$ evaluations,
we have $N_f^{(0)}=\mathcal O(n)$ and the total evaluation complexity is
$\mathcal O(n^2\varepsilon^{-2})$.
For a sample set of size $\bar m=\mathcal O(n^2)$, the bound becomes
$\mathcal O(n^3\varepsilon^{-2})$.

\subsection{Extension to matrix envelopes}
\label{subsec:matrix-envelopes}

The scalar analysis also applies to diagonal envelopes and to envelopes
constructed in an orthogonal frame. Orthogonal changes of coordinates
preserve the spectral curvature bounds.

\begin{corollary}
\label{cor:diagonal-envelope-models}
For each $k$, replace the scalar envelope models by
\[
  \underline q_k(x)
  =
  c_k+\frac12(x-x_k^+)^\top A_k(x-x_k^+),
  \qquad
  \overline q_k(x)
  =
  c_k+\frac12(x-x_k^+)^\top B_k(x-x_k^+),
\]
where $A_k,B_k\in\mathbb R^{n\times n}$ are symmetric and satisfy
\[
  0\preceq A_k\preceq B_k,
  \qquad
  \beta_{\min}I\preceq B_k
\]
for every $k$, with $\|B_k\|\leq\beta_{\max}$ whenever $\Delta_{k+1}<\Delta_k$.
The constants $\beta_{\min}>0$ and $\beta_{\max}<\infty$ are independent of $k$.
Replace the scalar touching and sufficient-decrease conditions by
\[
  f(x_k)
  =
  c_k+\frac12p_k^\top B_kp_k,
  \qquad
  f(x_k^+)
  \leq
  f(x_k)-\frac{\eta}{2}p_k^\top B_kp_k.
\]
Under the remaining assumptions of the corresponding results, the Lyapunov
decrease, radius decay, first-order stationarity, and function-evaluation
complexity bounds remain valid with $\beta_k$ replaced by $\|B_k\|$ and with
the same constants.
\end{corollary}

\begin{proof}
Since $B_k\succeq\beta_{\min}I$, the sufficient-decrease condition gives
\[
  f(x_k)-f(x_k^+)
  \geq\frac{\eta}{2}p_k^\top B_kp_k
  \geq\frac{\eta\beta_{\min}}{2}\|p_k\|^2.
\]
This is the same decrease bound used in the scalar analysis, so the
Lyapunov decrease and radius decay results hold with the same constants.
For each sample point, let $s_i=x_{k,i}-x_k$. The upper-envelope inequality
and the touching condition give
\[
\begin{aligned}
  0\leq f(x_{k,i})-f(x_k)
  &\leq p_k^\top B_ks_i+\frac12s_i^\top B_ks_i\\
  &\leq\|B_k\|
  \left(\|p_k\|\|s_i\|+\frac12\|s_i\|^2\right).
\end{aligned}
\]
The lower bound follows from the incumbent property, and the final
inequality uses the definition of the spectral norm $\|B_k\|$.
Applying the same weighting argument as in
Lemma~\ref{lem:envelope-rigidity} therefore gives
\[
  w_{k,i}|f(x_{k,i})-f(x_k)|
  \leq\frac52\|B_k\|\Delta_k^2.
\]
Thus the envelope estimates used in the first-order stationarity and
complexity proofs hold with $\beta_k$ replaced by $\|B_k\|$. The geometry condition and its
evaluation cost do not depend on the envelope matrices, so the remaining
results follow with the same constants.
\end{proof}

For an orthogonal matrix $Q_k$ satisfying $Q_k^\top Q_k=I$, let
\[
A_k=Q_k\operatorname{diag}(\alpha_k)Q_k^\top,
\qquad
B_k=Q_k\operatorname{diag}(\beta_k)Q_k^\top.
\]
If $0\leq\alpha_{k,j}\leq\beta_{k,j}$ and $\beta_{k,j}\geq\beta_{\min}$ for all $k,j$,
and $\max_j\beta_{k,j}\leq\beta_{\max}$ whenever $\Delta_{k+1}<\Delta_k$,
the corollary applies to any sequence of orthogonal frames.
The convergence and evaluation-complexity guarantees do not require
the frames to represent the exact curvature directions of the objective.

%% file: sections/numerical-experiments.tex
\section{Numerical experiments}
\label{sec:numerical-experiments}

This section presents the numerical evaluation of \qules{} and its variant
\qules{}-R, which adapts the orthogonal frame used to construct the envelopes.
We first describe the implementation, then introduce the benchmark problems
and comparison solvers, and finally assess solver performance using data
and performance profiles.

\subsection{Implementation}
\label{subsec:implementation}

The numerical experiments use a Julia implementation of
\qules{}\footnote{Code is available at
\url{https://github.com/open-dfo/QULES.jl}.}.
We formulate the inner problem in
JuMP~\cite{Lubin2023} and solve it locally with the nonlinear programming solver
Ipopt~\cite{WachterBiegler2006}. The inner problem is nonconvex, but the
analysis in Sect.~\ref{sec:theory} does not require a globally or locally
optimal solution.

Proposition~\ref{prop:inner-feasibility} gives a feasible start for the
mathematical subproblem. In the numerical implementation, we consider zero,
centroid, and sample-based trial displacements and keep those that admit scalar
curvatures satisfying the sampled envelope constraints. We pass the candidate
with the smallest curvature gap to Ipopt. The zero displacement admits
feasible curvatures at the incumbent in exact arithmetic.
If the inner solve fails, the implementation shrinks the radius and retries.
The convergence and complexity results concern Algorithm~\ref{alg:qules},
which requires a feasible pair with a gap no larger than the starting gap
at every iteration.

The numerical implementation uses curvature vectors,
\[
  \alpha_k=(\alpha_{k,1},\ldots,\alpha_{k,n})^\top,
  \qquad
  \beta_k=(\beta_{k,1},\ldots,\beta_{k,n})^\top,
\]
with $0\leq\alpha_{k,j}\leq\beta_{k,j}$. The corresponding envelope models
therefore use diagonal curvature matrices. The scalar formulation can be
constructed from $n+3$ sample points and represents isotropic curvature,
with spherical envelope level sets. The diagonal implementation starts with $2n+1$ sample
points and allows the curvature to vary across coordinate directions.
The initial sampling cost therefore grows linearly with $n$, consistent
with $N_f^{(0)}=\mathcal O(n)$.
Before solving the inner problem, the implementation scales the coordinates
and objective values and normalizes the trust-region radius to one to improve
numerical conditioning.
The inner problem then minimizes the sum of the diagonal curvature gaps in
the scaled coordinates.

\begin{figure}[t]
\centering
\includegraphics[width=\linewidth]{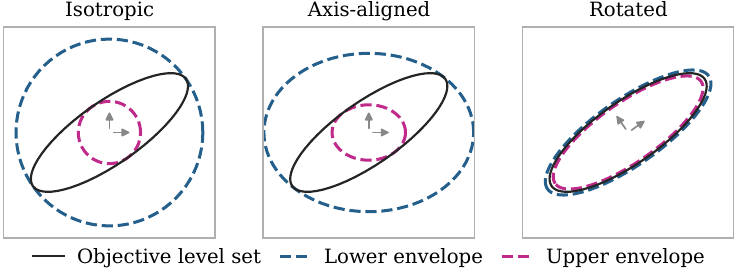}
\caption{Level sets of a rotated quadratic objective and its isotropic,
axis-aligned, and rotated quadratic envelopes. Aligning the envelope
frame with the principal axes produces a tighter pair in this example.}
\label{fig:directional-envelope-geometry}
\end{figure}

\qules{}-R constructs diagonal envelopes in a rotated frame and maps them
back to the original variables. This allows their
principal directions to adapt to the local curvature of the objective
rather than remain aligned with the coordinate axes.
Figure~\ref{fig:directional-envelope-geometry} illustrates the isotropic,
coordinate-aligned, and rotated envelope families.
The derivative-free frame-update rule is described in
Appendix~\ref{sec:frame-update}.
The frame update uses a target sample size of $n(n+3)/2$ to capture curvature
information across different directions, with $\bar m=\mathcal O(n^2)$.
However, \qules{}-R starts with the same $2n+1$ sample points as the diagonal
variant and collects additional samples during subsequent iterations.
The initialization count therefore remains $N_f^{(0)}=\mathcal O(n)$.
Geometry improvement requires at most $n$ additional evaluations per
iteration, and the total evaluation complexity bound is
$\mathcal O(n^3\varepsilon^{-2})$ under the dimension-independent bounds
stated in Sect.~\ref{sec:theory}.

The initial sample set consists of the starting point $x_0$ and the $2n$
coordinate points $x_0\pm\Delta_0e_j$, $j=1,\ldots,n$. This provides
directional information along each coordinate axis for the initial
diagonal envelope construction. Every subsequently evaluated point is
added to the sample set. Distant non-incumbent points are removed
to limit storage and keep the sampled data local to the current iterate.
The implementation also checks the weighted geometry condition and adds
at most $n$ sample points when it is not satisfied.

A run terminates when the objective-evaluation budget is exhausted or when
the trust-region radius falls below a prescribed tolerance. The algorithm
also stops when the objective values over the retained sample points
become numerically flat.

\subsection{Benchmark settings}
\label{subsec:test-problems}

We assess the evaluation efficiency, accuracy, and robustness of
\qules{} and \qules{}-R on a benchmark comprising 125 unconstrained
smooth optimization problem instances with different structural properties.
Let $\mathcal P$ denote the set of benchmark instances and $\mathcal S$
the set of solvers. For each problem $p\in\mathcal P$, let $f_p$ denote
its objective function and $n_p$ its dimension. The benchmark problems
come from the collections of Andrei~\cite{andrei2008unconstrained},
COCO/BBOB~\cite{hansen2021coco}, and Mor\'{e}, Garbow, and
Hillstrom (MGH)~\cite{more1981testing}.
An additional benchmark comprising 135 CUTEst problems is presented in
Appendix~\ref{sec:cutest}.

We compare \qules{} and \qules{}-R with three derivative-free optimization
solvers based on different algorithmic approaches.
\begin{enumerate}
\item \textbf{\texttt{NEWUOA}}~\cite{powell2006newuoa} is Powell's
model-based trust-region method using quadratic interpolation.
We use the implementation provided by PRIMA.jl~\cite{prima2023}.
\item \textbf{\texttt{NOMAD}}~\cite{Audet2022NOMAD4} implements mesh
adaptive direct search (MADS). We use the implementation provided by NOMAD.jl.
\item \textbf{\texttt{CMA-ES}}~\cite{Hansen2016Tutorial} is a randomized,
population-based method that adapts the covariance of its search
distribution. The runs use CMAEvolutionStrategy.jl.
\end{enumerate}
All solvers use the default settings of their implementations, with the
initial trust-region radii of \qules{}, \qules{}-R, and \texttt{NEWUOA},
the initial poll step of \texttt{NOMAD}, and the initial step size of
\texttt{CMA-ES} set to one. The stopping radius is $10^{-8}$ for
\qules{}, \qules{}-R, and \texttt{NEWUOA}. \texttt{NOMAD} and \texttt{CMA-ES}
retain their default stopping tolerances.

For problem $p$, we give each solver a budget of $50(n_p+1)$ objective
evaluations. We run each solver once except \texttt{CMA-ES}.
For \texttt{CMA-ES}, we use the pointwise median of the best-so-far
histories over ten fixed random seeds.

To define a common initialization procedure for the first three benchmark suites,
let $x_p^{\mathrm{orig}}$ denote the original reference minimizer of
problem $p$. For every problem in these suites, define the shifted reference minimizer
$x_p^\star$ componentwise by
\[
  (x_p^\star)_i=(x_p^{\mathrm{orig}})_i+0.5\sin(i),
  \qquad i=1,\ldots,n_p.
\]
Then define the starting point by
\[
  x_{p,0}=x_p^\star+\delta_p d_p,
\]
where $d_p$ is a predefined starting direction and
$\delta_p\in\{0.1,1,10,100\}$ whenever multiple starting-point scales are
used. The MGH problems use their standard starting points.

Solver performance is assessed using data profiles and performance
profiles~\cite{more2009benchmarking}.
For the data and performance profiles, let $f_p^L$ denote the lowest
objective value observed for problem $p$ over the compared solver histories.
A problem $p$ is considered solved by solver $s\in\mathcal S$ at
tolerance $\tau\in(0,1]$ once an evaluated point $x$ satisfies
\begin{equation}
  \frac{f_p(x)-f_p^L}{f_p(x_{p,0})-f_p^L}\leq\tau.
  \label{eq:benchmark-target}
\end{equation}
We use the convergence tolerances
$\tau\in\{10^{-3},10^{-7}\}$.

Data profiles measure the fraction of problems solved as a function of
the evaluation budget. Let $\nu$ denote the normalized evaluation budget. The data
profile is defined as
\begin{equation}
  d_s(\nu;\tau)=\frac{1}{|\mathcal P|}
  \left|\left\{p\in\mathcal P\mid
  \frac{t_{p,s}(\tau)}{n_p+1}\leq\nu\right\}\right|,
  \label{eq:benchmark-data-profile}
\end{equation}
where $t_{p,s}(\tau)$ is the first evaluation at which solver $s$
satisfies~\eqref{eq:benchmark-target} on problem $p$, and
$t_{p,s}(\tau)=+\infty$ if the tolerance is not reached within the
evaluation budget.

\begin{figure}[t]
\centering
\begingroup
\small
\captionof{table}{Structured benchmark suites. All solvers use the same starts and evaluation budgets.}
\label{tab:problem-groups}
\begin{tabular}{lrll}
\toprule
Group & Instances & $n$ & Curvature structure\\
\midrule
Well-conditioned separable & 40 & 2, 8, 15 & aligned curvature\\
Ill-conditioned separable & 25 & 2, 8, 15 & aligned, ratios up to $10^6$\\
Rotated ill-conditioned & 25 & 2, 8, 15 & same spectra, unaligned curvature\\
\bottomrule
\end{tabular}
\par
\endgroup
\medskip
\includegraphics[width=\linewidth]{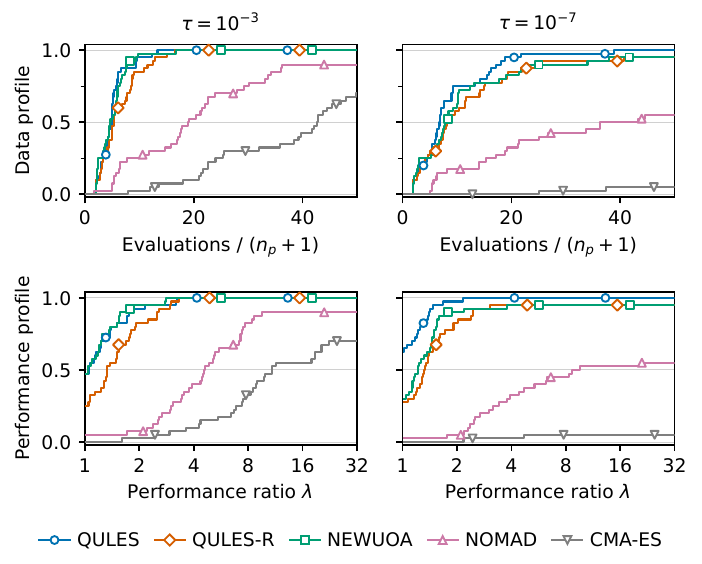}
\caption{40 Well-conditioned separable problems. Data profiles (top) 
and performance profiles (bottom) at $\tau=10^{-3}$ (left)
and $\tau=10^{-7}$ (right).}
\label{fig:separable-profiles}
\end{figure}
\begin{figure}[t]
\centering
\includegraphics[width=\linewidth]{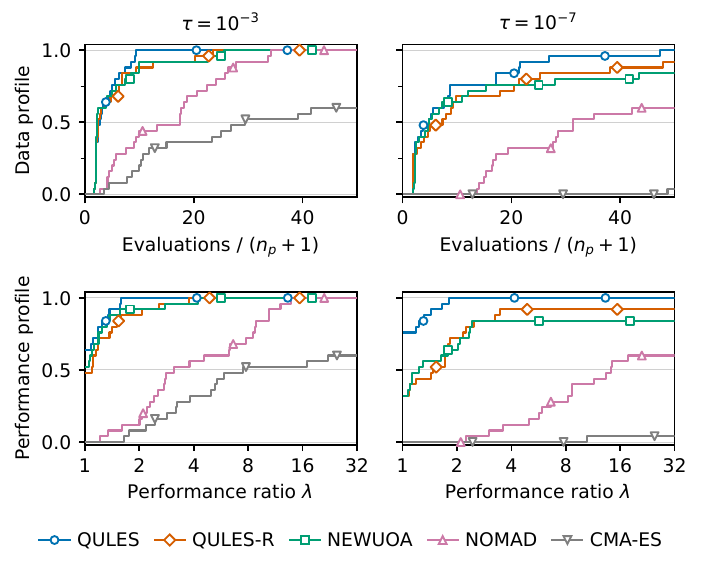}
\caption{25 Ill-conditioned separable problems. Data profiles (top) and performance profiles
(bottom) at $\tau=10^{-3}$ (left) and $\tau=10^{-7}$ (right).}
\label{fig:ill-conditioned-profiles}
\end{figure}
\begin{figure}[t]
\centering
\includegraphics[width=\linewidth]{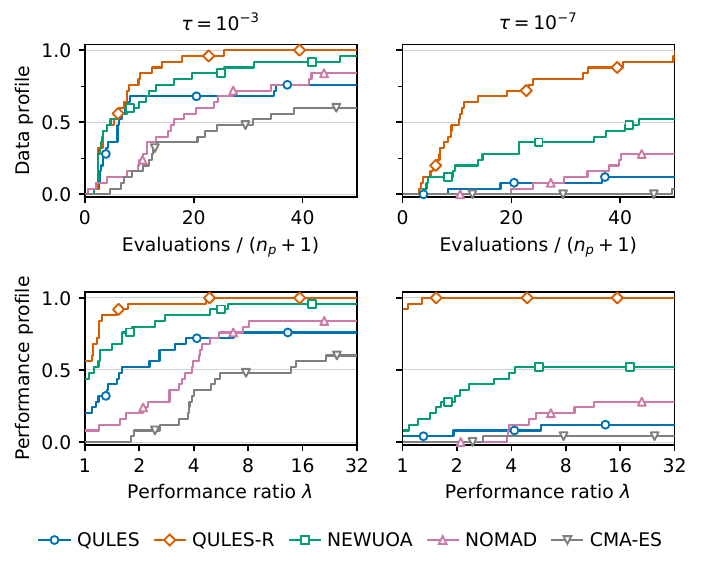}
\caption{25 Rotated ill-conditioned problems. Data profiles (top) and performance
profiles (bottom) at $\tau=10^{-3}$ (left) and $\tau=10^{-7}$ (right).}
\label{fig:rotated-profiles}
\end{figure}
Performance profiles compare solver efficiency relative to the best
solver on each problem. Define
\begin{equation}
  r_{p,s}(\tau)=
  \frac{t_{p,s}(\tau)}{\min_{s'\in\mathcal S}t_{p,s'}(\tau)}.
  \label{eq:benchmark-performance-ratio}
\end{equation}
If no solver satisfies the prescribed tolerance on problem $p$, we set
$r_{p,s}(\tau)=+\infty$ for every $s\in\mathcal S$. The performance
profile is
\begin{equation}
  \pi_s(\lambda;\tau)=\frac{1}{|\mathcal P|}
  \left|\left\{p\in\mathcal P\mid r_{p,s}(\tau)\leq\lambda\right\}\right|,
  \qquad \lambda\geq1.
  \label{eq:benchmark-performance-profile}
\end{equation}
Thus $\pi_s(\lambda;\tau)$ gives the fraction of problems for which
solver $s$ is within a factor $\lambda$ of the best solver at
tolerance $\tau$.

\subsection{Results}
\label{subsec:results}

The four benchmark suites cover well-conditioned separable, ill-conditioned separable,
rotated ill-conditioned, and general nonlinear problems (MGH)
(Tables~\ref{tab:problem-groups} and~\ref{tab:mgh-dimensions}).
We present the numerical results for each suite using data and performance
profiles.

\subsubsection{Well-conditioned separable problems}
Figure~\ref{fig:separable-profiles} shows that \texttt{NEWUOA}, \qules{}, and
\qules{}-R exhibit similar performance at low accuracy, reaching the target
within approximately $20(n_p+1)$ objective evaluations. At high accuracy,
methods exploiting curvature information show a clearer advantage over
direct-search methods. \qules{} performs slightly better than the other
solvers at $\tau=10^{-7}$, reaching all 40 targets, whereas
\texttt{NEWUOA} and \qules{}-R each reach 38. \qules{} also requires the fewest
evaluations on 25 of the 40 instances, compared with 12 for \texttt{NEWUOA},
including ties.
This performance is consistent with the diagonal envelope structure of
\qules{}, which is well suited to separable problems whose curvature is
aligned with the coordinate axes. \qules{}-R provides no clear additional
benefit in this setting.

\subsubsection{Ill-conditioned separable problems}
Figure~\ref{fig:ill-conditioned-profiles} shows that \qules{} and \qules{}-R
remain effective on separable problems with curvature ratios as large as
$10^6$. At low accuracy, \qules{}, \qules{}-R, and \texttt{NEWUOA} perform
similarly, each reaching all 25 targets and requiring fewer evaluations
than \texttt{NOMAD}. At high accuracy, \qules{} achieves the strongest results,
reaching all 25 targets compared with 21 for \texttt{NEWUOA} and 23
for \qules{}-R. \qules{} also requires the fewest evaluations on 19
instances, compared with 8 for \texttt{NEWUOA}, including ties. These results demonstrate
the robustness of the diagonal envelope construction to large differences
in curvature across coordinate directions. The independent curvature
parameters allow \qules{} to accommodate severe ill-conditioning while
maintaining an advantage over \texttt{NEWUOA}. \qules{}-R also remains
competitive, although adapting the coordinate frame offers no clear
additional benefit on these problems either.

\subsubsection{Rotated ill-conditioned problems}
The third benchmark suite is obtained by rotating the ill-conditioned
separable problems from the second suite. This introduces cross-coordinate
terms and removes the alignment between the curvature directions and the
coordinate axes. Figure~\ref{fig:rotated-profiles} shows a substantial decline
in the performance of \qules{} under this change. Although \qules{} initially
reduces the objective rapidly, zig-zagging slows its local convergence.
\texttt{NOMAD} therefore performs better than \qules{} on
these problems. \texttt{NEWUOA} performs well at low accuracy because its
quadratic interpolation model captures the cross-coordinate curvature.

\qules{}-R achieves the best performance on this suite, both in evaluation
efficiency and in the number of targets reached. At high accuracy,
\qules{}-R reaches all 25 targets, compared with 13 for
\texttt{NEWUOA} and 3 for \qules{}, and reaches the target in fewer
evaluations than every other solver on 22 instances. These results indicate that the diagonal envelope structure
itself is not the main limitation. \qules{}-R learns a frame adapted to the
curvature directions and maintains fast progress during local convergence.

\begin{figure}[t]
\centering
\begingroup
\small
\captionof{table}{Dimensions of the 35 MGH problems. Standard starting points are used.}
\label{tab:mgh-dimensions}
\begin{tabular}{lrrrrrrrrr}
\toprule
Dimension $n$ & 2 & 3 & 4 & 5 & 6 & 8 & 10 & 11 & 12\\
\midrule
Number of problems & 6 & 6 & 6 & 2 & 2 & 1 & 10 & 1 & 1\\
\bottomrule
\end{tabular}
\par
\endgroup
\medskip
\includegraphics[width=\linewidth]{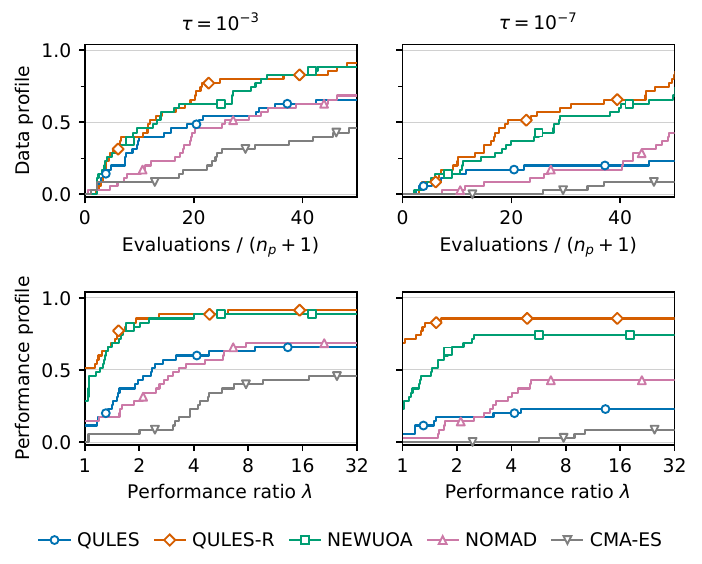}
\caption{35 MGH problems. Data profiles (top) and performance profiles (bottom)
at $\tau=10^{-3}$ (left) and $\tau=10^{-7}$ (right).}
\label{fig:mgh-profiles}
\end{figure}

\subsubsection{General nonlinear problems}
The MGH collection contains a range of smooth nonlinear problems with
different curvature structures, including nonseparable couplings, curved
valleys, and badly scaled variables. Figure~\ref{fig:mgh-profiles} shows that
\qules{}-R and \texttt{NEWUOA} perform well at low
accuracy, reaching 32 and 31 of the 35 targets, respectively.
\texttt{NOMAD} reaches 24 targets, while \texttt{CMA-ES} reaches 16.
At high accuracy, \texttt{NOMAD} reaches 15
targets, and \texttt{CMA-ES} reaches only 3 within the evaluation budget.

\qules{}-R achieves the strongest high-accuracy results, reaching 30
targets compared with 26 for \texttt{NEWUOA}, and requires the fewest
evaluations on 24 instances, compared with 8 for \texttt{NEWUOA}, including ties.
\qules{} reaches only 8 high-accuracy targets. This shows the limitation
of its coordinate-aligned frame. The frame-learning approach in \qules{}-R
addresses this limitation and achieves better performance than
\texttt{NEWUOA} without constructing a quadratic interpolation model.
The frame updates use previously evaluated points and require no
additional objective evaluations.

%% file: sections/conclusion.tex
\section{Conclusions}
\label{sec:conclusion}

We introduced \qules{}, a derivative-free trust-region algorithm that
generates trial points by solving a nonconvex envelope subproblem. This
subproblem constructs lower and upper quadratic envelopes consistent with
sampled objective values and minimizes their curvature gap, with their
common center defining the trial point. The convergence analysis requires
only a feasible solution of the subproblem, rather than a locally or
globally optimal one. We also introduced \qules{}-R, which adapts the
orthogonal envelope frame from previously evaluated points.

The trust-region radii converge to zero whenever the objective values along
the iterates are bounded below. For continuously differentiable objectives
whose initial lower level set is compact,
we show that $\|\nabla f(x_k)\|\to 0$ for any sequence of orthogonal frames.
For this result, it suffices to bound the upper-envelope curvature uniformly
on iterations that reduce the trust-region radius. Under the same curvature
bound and with a Lipschitz-continuous gradient, we obtain an
$\mathcal{O}(n^2\varepsilon^{-2})$ evaluation complexity bound for
\qules{} that includes geometry repair, provided that the sample set
size is $\mathcal O(n)$ and the remaining constants are independent
of $n$. For \qules{}-R, which retains $\mathcal O(n^2)$ samples, the
complexity bound becomes $\mathcal{O}(n^3\varepsilon^{-2})$.

Numerical experiments demonstrate the effectiveness of \qules{} on
coordinate-aligned problems. Within $50(n+1)$ objective evaluations per
problem, \qules{} and \texttt{NEWUOA} reach every target on both separable
suites at $\tau=10^{-3}$. At $\tau=10^{-7}$, \qules{} reaches all 40 targets
on the well-conditioned suite and all 25 on the ill-conditioned suite,
compared with 38 and 21 for \texttt{NEWUOA}. \qules{}-R reaches 38 and 23,
respectively, so adapting the frame offers no gain in target coverage on
these coordinate-aligned problems.

On the rotated ill-conditioned problems, fixed-frame \qules{} reaches
only 3 of the 25 high-accuracy targets. \qules{}-R, by contrast, reaches all
25, compared with 13 for \texttt{NEWUOA}, and reaches the target in fewer
evaluations than every other solver on 22 instances. \qules{}-R also reaches the most
targets on the 35 general nonlinear MGH problems. At $\tau=10^{-3}$, it
reaches 32 targets, compared with 31 for \texttt{NEWUOA}, 24 for
\texttt{NOMAD}, and 16 for \texttt{CMA-ES}. At $\tau=10^{-7}$, \qules{}-R
reaches 30 targets, followed by \texttt{NEWUOA} with 26, \texttt{NOMAD} with
15, fixed-frame \qules{} with 8, and \texttt{CMA-ES} with 3.

The present work focuses on unconstrained optimization. The formulation
allows us to adapt the envelope subproblem to
different problem structures. Future work could investigate how to use
this flexibility to incorporate bound, linear, and nonlinear constraints.

%% file: sections/frame-update.tex
\section{Derivative-free frame updates}
\label{sec:frame-update}

This section describes the frame update used by \qules{}-R. The update acts
on the orthogonal frame, independently of the envelope curvature vectors.
Consider displacements $d_\ell$ from the current incumbent, including the
newly evaluated trial, and define
\[
 D=[d_1,\ldots,d_m],\qquad
 \delta_\ell=f(x_k+d_\ell)-f(x_k).
\]
For an exactly quadratic objective with gradient $g$ and Hessian $H$ at
$x_k$,
\[
 \delta_\ell=g^\top d_\ell+\tfrac12d_\ell^\top H d_\ell.
\]
A weight vector $w$ satisfying $Dw=0$ eliminates $g$:
\begin{equation}
 2w^\top\delta
 =\left\langle H,D\operatorname{diag}(w)D^\top\right\rangle_F.
 \label{eq:supp-weighted-curvature}
\end{equation}
Thus a specified curvature functional can be recovered without constructing
all entries of $H$ if its target matrix lies in the range of the quadratic
moment map restricted to $\ker D$.

For columns $q_i,q_j$ of the current frame $Q$, use the targets
\[
 T_{ij}^{\mathrm{off}}=q_iq_j^\top+q_jq_i^\top,
 \qquad
 T_{ij}^{\mathrm{diag}}=q_iq_i^\top-q_jq_j^\top.
\]
Let $N$ have orthonormal columns spanning $\ker D$. For each target, the
implementation solves
\[
 \min_z\left\|D\operatorname{diag}(Nz)D^\top-T\right\|_F
\]
by an SVD pseudoinverse, choosing the minimum-norm solution, and sets $w=Nz$.
When both targets are exactly
representable, \eqref{eq:supp-weighted-curvature} gives
\[
 u_{ij}=2(w_{ij}^{\mathrm{off}})^\top\delta=2q_i^\top Hq_j,
 \qquad
 v_{ij}=2(w_{ij}^{\mathrm{diag}})^\top\delta
       =q_i^\top Hq_i-q_j^\top Hq_j.
\]
The angle
\[
 \vartheta_{ij}=\tfrac12\operatorname{atan2}(u_{ij},v_{ij})
\]
defines a Givens rotation $G_{ij}$ whose $(i,j)$ block is
\[
 \begin{bmatrix}
 \cos\vartheta_{ij}&-\sin\vartheta_{ij}\\
 \sin\vartheta_{ij}& \cos\vartheta_{ij}
 \end{bmatrix}.
\]
A sweep processes the pairs
$(1,2),(1,3),\ldots,(n-1,n)$ in lexicographic order, using the same sampled
differences throughout. At each pair, the sweep updates the frame to
$QG_{ij}$. Each rotation preserves orthogonality.

%% file: sections/cutest-appendix.tex
\section{CUTEst problems}
\label{sec:cutest}

We test the same five solvers on 135 unconstrained
CUTEst problems~\cite{gould2015cutest} with $2\le n\le 12$, using their
standard starting points (Table~\ref{tab:cutest-dimensions}) and the
evaluation budget of $50(n+1)$ and profile definitions
of Sect.~\ref{sec:numerical-experiments}.
Figure~\ref{fig:cutest-profiles} shows the resulting data and performance profiles.